\documentclass[11pt]{article}
\usepackage[all]{xy}
\usepackage{amsfonts,amsmath,oldgerm,amssymb,amscd}
\newcommand{\ra}{\rightarrow}		
\newcommand{\by}[1]{\stackrel{#1}{\ra}}

\newcommand{\ol}{\overline}		\newcommand{\wt}{\widetilde}
\newcommand{\iso}{\by \sim}

\newtheorem{theorem}{Theorem}[section]
\newtheorem{proposition}[theorem]{Proposition}
\newtheorem{lemma}[theorem]{Lemma}

\newtheorem{corollary}[theorem]{Corollary}
\newtheorem{question}[theorem]{Question}

\newcommand{\ga}{\alpha}	
		
\newcommand{\gf}{\varphi}

\newcommand{\gt}{\theta}

\newcommand{\gL}{\Lambda}

\newcommand{\BC}{\mbox{$\mathbb C$}}	
	\newcommand{\BF}{\mbox{$\mathbb F$}}

	\newcommand{\BR}{\mbox{$\mathbb R$}}

\newcommand{\CC}{\mbox{$\mathcal C$}}

	\newcommand{\CJ}{\mbox{$\mathcal J$}}

\newcommand{\mm}{\mbox{$\mathfrak m$}}	
	\newcommand{\p}{\mbox{$\mathfrak p$}}
\newcommand{\mq}{\mbox{$\mathfrak q$}}

\newcommand{\ot}{\mbox{\,$\otimes$\,}}	\newcommand{\op}{\,\mbox{$\oplus$\,}}
\newcommand{\Spec}{\mbox{\rm Spec\,}}	\newcommand{\hh}{\mbox{\rm ht\,}}
	
\newcommand{\Aut}{\mbox{\rm Aut}}	
\newcommand{\Hom}{\mbox{\rm Hom}}	\newcommand{\Sing}{\mbox{\rm Sing\,}}

\newcommand{\Um}{\mbox{\rm Um}}		\newcommand{\SL}{\mbox{\rm SL}}
		
\newcommand{\tr}{\mbox{\small \BR}}

\begin{document}

\begin{center}

{\LARGE Stability result for projective modules 
over blowup rings}
\vspace{.1in}

{\large Manoj Kumar Keshari} 

\vspace{.1in}

{\small \it School of Mathematics, Harish-Chandra Research Institute, \\
Chhatnag Road, Jhusi, Allahabad - 211019, India. \\
E-mail :  keshari@mri.ernet.in }
\end{center}

\section{Introduction}

Let $R$ be a normal affine domain of dimension $n\geq 3$
over an algebraically closed field $k$. Suppose {\rm char} $k=0$ or
{\rm char} $k=p \geq n$. Let
$g,f_1,\ldots,f_r$ be a $R$-regular sequence and  
$A=R[f_1/g,\ldots,f_r/g]$. Let $P$ be a stably free  $A$-module of
rank $n-1$. Then, Murthy proved that
there exists a projective $R$-module $Q$ such that $Q\ot_R A =P$ and
$\wedge^{n-1} Q=R$ (\cite{Mu}, Theorem 2.10).
As a consequence of Murthy's result, if $f,g\in \BC[X_1,\ldots,X_n]$
with $g\neq 0$, then all stably free 
modules over $\BC[X_1,\ldots,X_n,f/g]$ of rank $\geq n-1$ 
are free (\cite{Mu}, Corollary 2.11).

In this paper, we prove the following result (\ref{main2}),
which generalizes the above result of Murthy.
 
\begin{theorem}\label{1a}
Let $R$ be an affine algebra of dimension $n\geq 3$
over an algebraically closed field $k$. Suppose {\rm char} $k =0$ or
{\rm char} $k =p\geq n$. Let
$g,f_1,\ldots,f_r$ be a $R$-regular sequence and
$A=R[f_1/g,\ldots,f_r/g]$. Let $P'$ be a projective $A$-module of rank
$n-1$ which is extended from $R$. 
Let $(a,p) \in \Um(A \op P')$ and $P=A\op
P'/(a,p)A$.  Then, $P$ is extended from $R$.
\end{theorem}

Now, we will describe our next result. 
Let $R$ be an affine algebra over $\BR$ of dimension $n$. Let
$g\in R$ be an element not belonging to any real maximal ideal.  Let
$Q$ be a projective $R$-module of rank $\geq n-1$. Let $(a,p)\in
\Um(R_g\op
Q_g)$ and $P=R_g\op Q_g/(a,p)R_g$. Then, $P$ is extended from $R$
(\cite{Kesa}, Theorem 3.10).  This result was 
proved earlier by Ojanguren and Parimala in case $Q$ is free
(\cite{O-P}, Theorem).

In this paper, we prove the following result (\ref{main4}),
which is similar to (\ref{1a}).

\begin{theorem}
Let $R$ be an affine algebra of dimension $n\geq 3$
over $\BR$. Let $g,f_1,\ldots,f_r$ be a $R$-regular sequence and
$A=R[f_1/g,\ldots,f_r/g]$. Assume that $g$ does not belong to any real
maximal ideal of $R$. Let $P'$ be a projective $A$-module of rank
$\geq n-1$ which is extended from $R$.  Let $(a,p) \in \Um(A \op P')$ 
and $P=A\op P'/(a,p)A$.  Then, $P$ is extended from $R$.

In particular, every stably free $A$-module of rank
$n$ is extended from $R$.
\end{theorem}

As a consequence of above result, if
$f,g\in \BR[X_1,\ldots,X_n]$ with $g$ not belonging to any real
maximal ideal, then all stably free  modules of rank
$\geq n-1$ over $\BR[X_1,\ldots,X_n,f/g]$ are free (\ref{a3}).

The proof of the main theorem makes crucial use of results and
techniques of \cite{Bh}.


\section{Preliminaries}

In this paper, all the rings are assumed to be commutative Noetherian
and all the projective modules are finitely generated. We denote the
Jacobson radical of $A$ by $\CJ(A)$.

Let $B$ be a ring and let $P$ be a projective $B$-module. Recall that $p\in
P$ is called a {\it unimodular element} if there exists a $\psi \in
P^*=\Hom_B(P,B)$ such that $\psi(p)=1$. We denote by $\Um(P)$, the set
of all unimodular elements of $P$. We write $O(p)$ for the ideal of
$B$ generated by $\psi(p)$, for all $\psi \in P^*$. Note that, if
$p\in \Um(P)$, then $O(p)=B$.
  
Given an element $\gf\in P^\ast$ and an element $p\in P$, we define an
endomorphism $\gf_p$ of $P$ as the composite $P\by \gf B\by p P$.  If
$\gf(p)=0$, then ${\gf_p}^2=0$ and, hence, $1+\gf_p$ is a unipotent
automorphism of $P$.

By a {\it transvection}, we mean an automorphism of $P$ of the form
$1+\gf_p$, where $\gf(p)=0$ and either $\gf \in \Um(P^\ast)$
or $p \in\Um(P)$. We denote by $E(P)$, the subgroup of
$\Aut(P)$ generated by all transvections of $P$. Note that, $E(P)$ is a
normal subgroup of $\Aut(P)$. 

An existence of a transvection of $P$ pre-supposes that $P$ has a
unimodular element.
Let $P = B\op Q$, $q\in Q, \alpha\in
Q^*$. Then, the automorphisms $\Delta_q$ and $\Gamma_\ga$ of $P$ defined
by $\Delta_q(b,q')=(b,q'+bq)$ and
$\Gamma_\alpha(b,q')=(b+\alpha(q'),q')$ are transvections of
$P$. Conversely, any transvection $\Theta$ of $P$ gives rise to a
decomposition  $P=B\op Q$ in such a way that $\Theta = \Delta_q$ or
$\Theta = \Gamma_\alpha$.

\begin{define}
Let $A$ be a ring and let $P$ be a projective $A$-module. We say that $P$
is {\it cancellative} if $P\op A^r \simeq Q\op A^r$ for some positive
integer $r$ and some projective $A$-module $Q$ implies that $P\simeq Q$. 
\end{define}

We begin by  stating two classical results due to Serre \cite{S1} and
Bass \cite{B} respectively.

\begin{theorem}\label{serre}
Let $A$ be a ring with $\dim A/\CJ(A) =d$. Then, any projective
$A$-module $P$ of rank $> d$ has a unimodular element. 
\end{theorem}   


\begin{theorem}\label{bass1}
Let $A$ be a ring of dimension $d$ and let $P$ be a projective
$A$-module of rank $>d$. Then $E(A\op P)$ acts transitively on
$\Um(A\op P)$. 
In particular, $P$ is cancellative. 
\end{theorem}

The above result of Bass is best possible in general. But, in case of
affine algebras over algebraically closed fields, we have the
following result  due to Suslin \cite{Su1}.

\begin{theorem}\label{007}
Let $A$ be an affine algebra of dimension $n$ over an algebraically closed
field. Then, all projective $A$-modules of rank $\geq n$ are cancellative. 
\end{theorem}

\begin{remark}\label{lift}
Let $P$ be a finitely generated
projective $A$-module of rank $d$. Let $t$ be a non zero
divisor of $A$ such that $P_t$ is free. Then, it is easy to see that
there exits a free submodule $F=A^d$ of $P$ and a positive integer $l$
such that, if $s=t^l$, then $sP\subset F$. Therefore, $sF^*\subset
P^*\subset F^*$.
If $p\in F$, then $\Delta_p\in E(A\op F)
\cap E(A\op P)$ and if $\alpha \in F^*$, then $\Gamma_{s\alpha}\in
E(A\op F)\cap E(A\op P)$. 
\end{remark}

\medskip

The following result is due to Bhatwadekar and Roy ($\cite{Bh-Roy}$,
Proposition 4.1).

\begin{proposition}\label{one}
Let $A$ be a ring and let $J$ be an ideal of $A$. Let $P$ be a
projective $A$-module of rank $n$. 
Then, any transvection $\wt \Theta$ of $P/JP$, i.e. $\wt \Theta \in E(P/JP)$,
can be lifted to a (unipotent) automorphism $\Theta$ of $P$. In
particular, if $P/JP$ is free of rank $n$, then any element $\ol \Psi$ of
$E((A/J)^n)$ can be lifted to $\Psi \in \Aut(P)$. If, in
addition, the natural map $\Um(P) \ra \Um(P/JP)$ is surjective, then
the natural map $E(P)\ra E(P/JP)$ is surjective.
\end{proposition}

\begin{define}
For a ring $A$, we say that projective stable range of $A$ is $\leq r$
(notation: psr$(A)\leq r$) if for all projective $A$-modules $P$ of rank
$\geq r$ and $(a,p)\in \Um(A\op P)$, we can find  $q\in P$ such that
$p+aq\in \Um(P)$. Similarly, $A$ has stable range
$\leq r$ (notation: sr$(A)\leq r$) is defined the same way as
psr$(A)$ but with $P$ required to be free. 
\end{define}

The following result is due to Bhatwadekar ($\cite{Bh}$, Corollary 3.3)
and is a generalization of a result of Suslin ($\cite{Su2}$, Lemma 2.1).
See \cite{Kesa}, for the definition of $ESp_4(B)$.

\begin{proposition}\label{two}
Let $B$ be a ring with ${\rm psr}(B) \leq 3$
and let $I$ be an ideal of $B$.
Let $P$ be a projective $B$-module of
rank $2$ such that $P/IP$ is free. Then, any element of
$\SL_2(B/I) \cap ESp_4(B/I)$ can be lifted to an element of $\SL(P)$.
\end{proposition}

\begin{remark}
In (\cite{Bh}, Corollary 3.3), Proposition \ref{two} is stated with the
assumption that $\dim B=2$. However, the proof works equally well in above
case.
\end{remark}

The following result is due to Mohan Kumar, Murthy and Roy (\cite{MMR},
Theorem 3.7) and is used in (\ref{mmr1}). 

\begin{theorem}\label{mmr}
Let $A$ be an affine algebra of dimension $d\geq 2$ over $\ol
\BF_p$. Suppose that $A$ is regular when $d=2$. Then ${\rm psr}(A) \leq d$.
\end{theorem}

The following two results are due to Lindel (\cite{L1}, Theorem) and
(\cite{L}, Theorem 2.6). Recall that a ring $A$ is called essentially
of finite type over a field $k$, if $A$ is the localisation of an
affine algebra over $k$.
     
\begin{theorem}\label{lindel1}
Let $A$ be a regular ring which is essentially of finite type over a
field $k$. Then, every projective $A[X]$-module is extended from $A$.
\end{theorem}
 
\begin{theorem}\label{lindel}
Let $B$ be a ring of dimension $d$ and let $R =
B[T_1,\ldots,T_n]$. Let $P$ be a projective $R$-module of rank $\geq$
max $(2,d+1)$. Then $E(P\op R)$ acts transitively on $\Um(P\op R)$.
\end{theorem}

The following result is due to Suslin (\cite{Su3}, Theorem 2). The special
case, namely $n=2$ was proved earlier by Swan and Towber \cite{S-T}.

\begin{theorem}\label{c1}
Let $A$ be a ring and $[a_0,a_1,\ldots,a_n]\in \Um_{n+1}(A)$.
 Then, there exists $\Gamma\in \SL_{n+1}(A)$ with
$[{a_0}^{n!},a_1,\ldots,a_n]$ as the first row.
\end{theorem}

The next three results are due to Suslin (\cite{Su2}, Proposition 1.4,
1.7 and Corollary 2.3) and are very crucial for the proof of our main
theorem (see also \cite{Mu}, Remark 2.2).
Here ``cd'' stands for cohomological dimension (See
\cite{Serre} for definition).

\begin{proposition}\label{c2}
Let $X$ be a regular affine curve over a field $k$ and let $l$ be a
prime with $l\neq$ {\rm char} $k$. Suppose that ${\rm cd}_l \,k \leq 1$.
Then, the group $SK_1(X)$ is $l$-divisible.
\end{proposition}

\begin{proposition}\label{c3}
Let $X$ be a regular affine curve over a field $k$ of characteristic
$\neq 2$ and ${\rm cd}_2 \,k\leq 1$. Then, the 
canonical homomorphism $K_1Sp(X) \ra SK_1(X)$ is an isomorphism.
\end{proposition}

\begin{proposition}\label{b1}
Let $A$ be a ring and $[a_1,\ldots,a_n] \in \Um_n(A)$
($n\geq 3$). Let $I=\sum_{i\geq 3} Aa_i$ and
$J=\sum_{i\geq 4} Aa_i$ be ideals of $A$. Let $b_1,b_2\in A$ be such that
$Ab_1+Ab_2+I=A$. Let ``bar'' denotes reduction mod $I$. Suppose that 

$(i)\; \dim A/I \leq 1$ and ${\rm sr}(A/J) \leq 3$,

$(ii)$ there exists an $\ol \ga\in \SL_2(\ol A)\cap ESp(\ol A)$, such
that $[\ol a_1,\ol a_2]\ol \ga = [\ol b_1,\ol b_2]$. \\
Then, there exists a $\gamma\in E_n(A)$ such that $[a_1,\ldots,a_n]\gamma = 
[b_1,b_2,a_3\ldots,a_n]$.
\end{proposition}

Using above results, Suslin proved the following cancellation
theorem (\cite{Su2} Theorem 2.4).

\begin{theorem}\label{cancel}
Let $A$ be an affine algebra of dimension $d\geq 2$ over an infinite
perfect field $k$. Suppose ${\rm cd} \,k \leq 1$ and $d! \in k^*$. Let
$[a_0,a_1,\ldots,a_d]\in \Um_{d+1}(A)$ and let $r$ be a
positive integer. Then, there exists $\Gamma \in E_{d+1}(A)$ such that
$[a_0,a_1,\ldots,a_d]\Gamma = [{c_0}^r,c_1,\ldots,c_d]$. As a
consequence, every stably free $A$-module of rank $d$ is free (Theorem
\ref{c1}).
\end{theorem}

The following result is due to Bhatwadekar (\cite{Bh}, Theorem 4.1)
and is a generalisation of above result of Suslin.

\begin{theorem}\label{b2}
Let $A$ be an affine algebra of dimension $d\geq 2$ over an infinite
perfect field $k$. Suppose ${\rm cd}\, k \leq 1$ and $d! \in k^*$.
Then, every projective $A$-module $P$ of rank $d$ is cancellative. 
\end{theorem}

The following result is very crucial for our main theorem and the
proof of it is contained in (\cite{Bh}, Theorem 4.1). 

\begin{proposition}\label{b5}
Let $A$ be a ring and let $P$ be a projective $A$-module
of rank $d$. Let $s\in A$ be a non-zero-divisor such that $P_s$ is
free. Let $F=A^d$ be a free submodule of $P$ with $F_s=P_s$ and
$sP\subset F$. Let $e_1,\ldots,e_d$ denote the standard basis of $F$.
Let $(a,p)\in \Um(A\op P)$ be such that 

$(1)~ a=1$ mod $As$,

$(2)~ p\in F\subset P$ with 
$p={c_1}^de_1+c_2e_2+\ldots+c_de_d$, for some $c_i\in A$,

$(3)$ every stably free $A/Aa$-module of rank $\geq
d-1$ is free. \\
Then, there
exists $\Delta\in \Aut(A\op P)$ such that $(a,p)\Delta=(1,0)$.
\end{proposition}

The following result is used to prove our second result
(\ref{main4}) and is due to Ojanguren and Parimala ($\cite{O-P}$,
Proposition 3 and Proposition 4). 

\begin{proposition}\label{op}
Let $\CC = \Spec C$ be a smooth affine curve over a field $k$ of
characteristic $0$. Suppose that every residue field of $\CC$ at a
closed point has cohomological dimension $\leq 1$. Then, $SK_1(C)$ is
divisible and the
natural homomorphism $K_1Sp(C) \ra SK_1(C)$ is an isomorphism..
\end{proposition}


\section{Main Theorem 1}

In this section, we will prove our first result (\ref{main2}). 
We begin with the following result, 
the proof of which is similar to (\cite{Mu}, Corollary 2.8).
 
\begin{lemma}\label{a1}
Let $R$ be an affine algebra of dimension $n\geq 3$
over a field $k$. Let
$g,f_1,\ldots,f_r$ be a $R$-regular sequence and
$A=R[f_1/g,\ldots,f_r/g]$. Let $P$ be a projective $A$-module of
rank $\geq n-1$ and $(a,p) \in \Um(A\op P)$. 
Then, there exists  $\Psi \in \SL(A\op P)$ such that 
$(a,p)\Psi = (1,0)$ mod $Ag$.
\end{lemma}

\begin{proof}
Since $g,f_1,\ldots,f_r$ is a $R$-regular sequence,
$A=R[X_1,\ldots,X_r]/I$, where $I=(gX_1-f_1,\ldots,gX_r-f_r)$. Let
``bar'' denote reduction modulo $Ag$. Then 
$\ol A=\ol R[X_1,\ldots,X_r]$, where $\ol
R=R/(g,f_1,\ldots,f_r)$.
Since $\dim \ol R \leq n-2$, by (\ref{lindel}), there
exists $\ol \Psi
\in E(\ol A \op \ol P)$ such that $(\ol a,\ol p)\ol \Psi = (1,0)$. By
(\ref{one}), we can lift $\ol \Psi$ to $\Psi \in \SL(A\op P)$.
Hence, we have $(a,p)\Psi = (1,0)$ mod $Ag$. 
$\hfill \square$
\end{proof}

\begin{lemma}\label{a11}
Let $R$ be an affine algebra of dimension $n\geq 3$ over a field
$k$. Let $g,f_1,\ldots,f_r \in R$ with $g$ a non-zero-divisor and 
$A=R[f_1/g,\ldots, f_r/g]$. Let $S=1+gR$ and $B=A_S$. Let $P$ be a
projective $B$-module of rank $\geq n-1$ which is extended from $R_S$. Let
$(a,p)\in \Um(B\op P)$ with $(a,p)=(1,0)$ mod $Bg$.  Then, we have the
followings:

(1) there exists $s\in R$ such that $P_s$ is free and

(2)  there exists $\Delta\in \Aut(B\op P)$ such that $(a,p)\Delta
=(1,0)$ mod $Bsg$.\\
Further, given any ideal $J$ of $R$ of height $\geq 1$ with
$\hh(g,J)R\geq 2$, 
we can choose $s$ such that $s\in J$.
\end{lemma}

\begin{proof}
Choose $s_1\in J$ such that $\hh s_1R=1$ and $\hh (s_1,g)R \geq 2$.
By replacing $(a,p)$ by $(a+\ga(p),p)$ for some $\ga\in P^*$, if
necessary, we may assume that $\hh aB=1$ and $\hh (a,s_1)B \geq 2$. 
Suppose $\p_1,\ldots,\p_t$ are minimal primes of
$gR_S$; $\wt \p_1,\ldots,\wt \p_{t'}$
 : minimal primes of $R_S$ and
$\mq_1,\ldots,\mq_{t''}$ : minimal primes of $aB$. Since $P$ is
extended from $R_S$,  $P_\Sigma$ is
free, where $\Sigma=R_S \backslash \cup_{i=1}^t\, \p_i
\cup_{j=1}^{t'}\, \wt \p_j \cup_{l=1}^{t''}\, (\mq_l\cap R_S)$ 
(any projective module over a semi local ring is free). 
Hence, there exists some $s_2\in \Sigma$ such
that $P_{s_2}$ is free. We may assume $s_2\in R$. Hence $\hh Rs_2 \geq
1$, $\hh (s_2,g)R_S \geq 2$ and $\hh (a,s_2)B\geq 2$.

Write $s=s_1s_2$. Then $\hh (g,s)R_S \geq 2$ and $\hh (a,s)B \geq  2$. 
Since $\dim R_S/(\CJ(R_S),s) \leq n-2$ 
and $P$ is extended from $R_S$, 
by (\ref{serre}), $P/sP$ has a unimodular
element. Write $B_1=B/Bs$, $P/sP=B_1\op P_1$ and $(b,p_1)$ as the
image of $p$ in $P/sP$.

Let ``bar'' denotes reduction modulo the ideal $B_1a$. Since $a=1$
mod $Bg$ and $\hh (s,a)B \geq 2$, 
$\dim \ol B_1\leq n-3$. Note that, $\ol p=(\ol b,\ol p_1) \in 
\Um(\ol P = \ol B_1 \op \ol P_1)$. 
Since rank of $\ol P\geq n-1$, by (\ref{bass1}),
there exists $\ol \Phi \in E(\ol P)$ such that $(\ol b,\ol p_1)\ol
\Phi=(1,0)$. 
In general, the natural map $\Um(B_1\op P_1) \ra 
\Um(\ol B_1\op \ol P_1)$ is surjective. Hence, by 
(\ref{one}), we can lift $\ol \Phi$ to some element $\Phi\in E(B_1\op P_1)$.
Let $(b,p_1)\Phi=(c,p_2)$. Then $(c,p_2) = (1,0)$ mod $B_1a$. 

Let $(c,p_2)=(1,0)-a(c_1,p_3)$ for some $(c_1,p_3) \in B_1\op P_1$. Then
$(a,c,p_2)\Delta_{(c_1,p_3)} =
(a,1,0)$, where $\Delta_{(c_1,p_3)} \in E(B_1\op P_1)$. 
Recall that $P/sP=B_1\op P_1$.
By (\ref{one}), we can lift $(1,\Phi) \Delta_{(c_1,p_3)} \in
E(B\op P/s(B\op P))$ to some element $\Psi \in \Aut (B\op P)$ such that
$(a,p)\Psi =(a,q)$ with $O(q)=B$ mod $Bs$. Since $a=1$ mod $Bg$, 
there exists $\Psi_1 \in E(B\op P)$ such that
$(a,q)\Psi_1 =(1,0)$ mod $Bsg$. Let $\Delta=\Psi\Psi_1$. Then
$(a,p)\Delta=(1,0)$ mod $Bsg$. This proves the lemma.
$\hfill \square$
\end{proof}


\begin{lemma}\label{main1}
Let $R$ be an affine algebra of dimension $n\geq 3$ over an
algebraically closed field $k$. Suppose {\rm char} $k=0$ or {\rm char}
$k=p \geq n$. Let $g,f_1,\ldots,f_r\in R$ with $g$ a non-zero-divisor
and $A=R[f_1/g,\ldots,f_r/g]$. Let $S=1+Rg$ and $B=A_S$.  Let $P$ be a
projective $B$-module of rank $n-1$ which is extended from $R_S$. Let
$(a,p)\in \Um(B\op P)$ with $(a,p)=(1,0)$ mod $Bg$. Then, there exists
$\Delta \in \Aut(B\op P)$ such that $(a,p)\Delta=(1,0)$.
\end{lemma}

\begin{proof}
Without loss of generality, we can assume that $R$ is reduced. 
Let $J_1$ be the ideal of $R_g$ defining the
singular locus $\Sing R_g$.
Since $R_g$ is reduced, $\hh J_1 \geq 1$. Note that $\sqrt {J_1}=J_1$.
Let $J=J_1\cap R$. Then, we may assume that $g$ does
not belong to any minimal primes of $J$ and $\hh J \geq 1$. Hence
$\hh(g,J)R\geq 2$.

Since $(a,p)=(1,0)$ mod $Bg$, by (\ref{a11}), there exists some
$s\in J$ and $\Phi\in \Aut(B\op P)$ such that $P_s$ is free and 
$(a,p)\Phi=(1,0)$ mod $Bsg$. Hence, replacing $(a,p)$ by $(a,p)\Phi$,
we can assume that $(a,p)=(1,0)$ mod $Bsg$. It is easy to see that
 we can replace $B$
by $C=A_T$, where $T=(1+gk[g])h$ for some $h\in 1+gR$. Note that, $B=C_S$.

Since $P_s$ is free of rank $n-1$, there exists a free
submodule $F=C^{n-1}$ of $P$ such that $F_s=P_s$. By replacing $s$ by
a power of $s$, we may assume that $sP\subset F$. Let
$e_1,\ldots,e_{n-1}$ denote the standard basis of $C^{n-1}$.
Since $(a,p)=(1,0)$ mod $Csg$, $p\in sgP\subset gF$. Let $p=b_1
e_1+\ldots+b_{n-1}e_{n-1}$, for some $b_i\in gC$. 
Then $[a,b_1,\ldots,b_{n-1}]\in \Um_n(C)$. As
$1-a \in Csg$, $[a,sb_1,\ldots,sb_{n-1}] \in \Um_n(C)$.\\

For $n=3$, by Swan's Bertini theorem (\cite{Swan}, Theorem 1.3) as quoted 
in (\cite{Mu1}, Theorem 2.3),
there exists $c_1,c_2 \in C$ such that, if  $a'=a+sb_1c_1+sb_2c_2$,
then $C/Ca'$ is a reduced regular (since $a'=1$ mod
$Csg$ and $s\in J$) affine $k(g)$-algebra of dimension $1$. 
For $n\geq 4$, by prime avoidance, 
there exists $c_1,\ldots,c_{n-1}\in C$ such that, if $a' =
a+sb_1c_1+\ldots+sb_{n-1}c_{n-1}$, then $C/Ca'$ is affine $k(g)$-algebra 
of dimension $\leq n-2$. Note that, $a'=1$ mod $Csg$.

Let ${e_1}^*,\ldots,{e_{n-1}}^*$ be a dual basis of $F^*$ and let $\gt_i =
c_i{e_i}^* \in F^*$. Then, by (\ref{lift}), $\Gamma_{s\gt_i} \in
E(C\op P)$ and $\Gamma_{s\gt_i}(a,p)=(a+sb_ic_i,p)$. Hence, it
follows that there exists $\Psi_1\in E(C\op P)$ such that
$\Psi_1(a,p)=(a',p)$. 

Let ``bar'' denote reduction modulo $Ca'$. Since
$Ca'+Cs=C$ and $P_s=F_s$ is free, the inclusion $F\subset P$ gives
rise to equality $\ol F=\ol P$. In particular, $\ol P$ is free of rank
$n-1 \geq 2$ with a basis $\ol e_1,\ldots,\ol e_{n-1}$ and $\ol p \in
\Um(\ol P)$. 
Recall that $\ol C$ is an affine algebra of dimension $n-2$ over a
$C_1$-field 
$k(g)$. Hence, by (\ref{b2}),
 every projective $\ol
C$-module of rank $n-2$ is cancellative.\\

If $n=3$, then $\ol C$ is a regular affine
algebra of dimension $1$ over a $C_1$-field $k(g)$. Hence, by
(\ref{c2}, \ref{c3}), $SK_1(\ol C)$ is a divisible group and
the canonical homomorphism $K_1Sp(\ol C) \ra SK_1(\ol C)$ is an
isomorphism. Hence, there exists $\Theta' \in
\SL_2(\ol C)\cap ESp(\ol C)$ and $t_1,t_2\in C$ such that, if $p_1 =
{t_1}^2e_1 +t_2e_2 \in F$, then $\Theta'(\ol p)=\ol p_1$. 
We have $\dim B/Bg=\dim A/Ag=2$ and $B_g$ is an $k(g)$-algebra of dimension
two. Thus $\Spec B=\Spec B/Bg \cup \Spec B_g$ with $\dim B/Bg=2=\dim
B_g$. Hence psr$(B)\leq 3$. Hence,
by (\ref{two}), $\Theta'\ot \ol B$ has a lift $\Theta\in \SL(P\ot B)$.

For $n\geq 4$. Since $\ol P$ is free of
rank $n-1$, $E_{n-1}(\ol C)=E(\ol P)$.
Hence, by (\ref{cancel}), there exists
$\wt \Theta \in E(\ol P)$ and $t_i\in C$, $1\leq i\leq n-1$
such that, if $p_1={t_1}^{n-1}e_1+t_2e_2+\ldots+t_{n-1}e_{n-1} \in F$,
then $\wt \Theta(\ol p) = \ol p_1$. By (\ref{one}), $\wt
\Theta$ can be lifted to an element $\Theta\in \SL(P)$.

Write $P$ for $P\ot B$. Thus, in either case, there exists $q \in P$ such that
$$\Theta(p)=p_1- a'q,~~ {\rm where}~~ p_1={t_1}^{n-1}e_1+t_2e_2+\ldots+
t_{n-1}e_{n-1}.$$ 

The automorphism
$\Theta$ of $P$ induces an automorphism $\gL_1=(Id_B,\Theta)$ of $B\op
P$. Let $\gL_2$ be the transvection $\Delta_{q}$ of $B\op P$. Then
$(a',p)\gL_1 \gL_2 =(a',p_1)$. 

By (\ref{b5}), there exists $\gL_3\in \Aut(B\op P)$ such that
$(a',p_1)\gL_3=(1,0)$. Let $\Delta=\Psi_1\gL_1\gL_2\gL_3$. Then $\Delta\in
\Aut(B\op P)$ and $(a,p)\Delta=(1,0)$. This proves the result.
 $\hfill \square$ 
\end{proof}

\begin{remark}
Let $A$ be a ring and  $g,h\in A$ with
$Ag+Ah=A$. Then, any projective $A$-module $E$ is given by a triple
$(Q,\ga,P)$, where $Q,P$ are projective modules over $A_h$ and $A_g$
respectively and $\ga$ is a prescribed $A_{gh}$-isomorphism $\ga:
Q_g\iso P_h$.

Let $g,h\in A$ with $Ag+Ah=A$ and let $P$ be a projective
$A$-module. Let $(a,p)\in \Um(A_g\op P_g)$ and 
$Q=A_g\op P_g/(a,p)A_g$. If $\gf : Q_h \iso P_{gh}$ is an isomorphism,
then the triple $(P_h,\gf,Q)$ yields a projective $A$-module $E$ such
that $Q=E\ot A_g$.
\end{remark}

Now, we prove the main result of this section. In case $P'$ is free
(i.e. $P$ is stably free), it is proved in (\cite{Mu}, Theorem 2.10).

\begin{theorem}\label{main2}
Let $R$ be an affine algebra of dimension $n\geq 3$
over an algebraically closed field $k$.
Suppose {\rm char} $k=0$ or {\rm char} $k=p \geq n$. Let
$g,f_1,\ldots,f_r$ be a $R$-regular sequence and
$A=R[f_1/g,\ldots,f_r/g]$. Let $P'$ be a projective $A$-module of rank
$n-1$ which is extended from $R$.
Let $(a,p)\in \Um(A\op P')$ and 
$P=A\op P'/(a,p)A$. 
Then, $P$ is extended from $R$.
\end{theorem}

\begin{proof}
By (\ref{a1}), there exists
$\Psi\in \SL(A\op P')$ such that $(a,p)\Psi = (1,0)$ mod $Ag$. Let $S=1+Rg$
and $B=A_S$. Applying (\ref{main1}), there exists $\Psi_1\in 
\Aut(B\op (P'\ot B))$ such that $(a,p)\Psi\Psi_1 = (1,0)$. Let $\Delta =
\Psi \Psi_1$. Then, there exists some $h\in 1+Rg$ such that $\Delta\in
\Aut(A_h \op {P'}_h)$ and $(a,p)\Delta=(1,0)$. We have the isomorphism
$\Gamma : P_h \iso {P'}_h$ induced from $\Delta$.
The module $P$ is given by the triple $({P'}_h, \Gamma_g,
 P_g)$. Since $Rg+Rh=R$, $R_g=A_g$, $R_{gh}=A_{gh}$ and $\Gamma_g :
P_{gh} \iso {P'}_{gh}$ is an isomorphism of $R_{gh}$ module, the
triple $({P'}_h,\Gamma_g,P_g)$ defines a projective $R$-module $Q$
of rank $n-1$ such that $P = Q\ot A$. This 
proves the theorem.
$\hfill \square$
\end{proof}
\medskip

The following result is a generalisations of (\cite{Mu}, Theorem
2.12), where it is proved for stably free modules. 

\begin{theorem}\label{mmr1}
Let $R$ be an affine domain of dimension $n\geq 4$
over $\ol \BF_p$. Suppose $p\geq n$. Let $K$ be the field of fractions
of $R$ and let $A$ be a sub-ring of $K$ with $R\subset A\subset K$.
Let $P'$ be a projective $A$-module of rank $n-1$ which is extended
from $R$. Let $(a,p)\in \Um(A\op P')$
and $P=A\op P'/(a,p)A$. Then, $P$ is extended from $R$.
\end{theorem}

\begin{proof}
We may assume that $A$ is finitely generated over $R$,
i.e. there exist $g,f_1,\ldots,f_r \in R$ such that 
$A=R[f_1/g,\ldots,f_r/g]$. 
Since $P'$ is extended from $R$, we can choose an element $s\in R$
such that ${P'}_s$ is free. Let ``bar'' denote reduction modulo $Asg$. Then 
$\ol A=A/Asg$ is an affine 
algebra of dimension $\leq n-1$ over $\ol \BF_p$. Since $n-1\geq 3$, 
by (\ref{mmr}),
psr$(\ol A) \leq n-1$. Hence, there exists $\ol \Psi \in E(\ol A\op \ol
P')$ such that $(\ol a,\ol p)\ol \Psi =(1,0)$. By 
(\ref{one}), $\ol \Psi$ can be lifted to $\Psi\in \SL(A\op
P')$. Replacing $(a,p)$ by $(a,p)\Psi$, we can assume that
$(a,p)=(1,0)$ mod $Asg$. Let $B=A_{1+gR}$. Then, by (\ref{main1})
there exists $\Gamma\in  \Aut(B\op (P'\ot B))$ such that
$(a,p)\Gamma=(1,0)$. 
Rest of the argument is same as in (\ref{main2}).
$\hfill \square$
\end{proof}

The following result is a generalisations of (\cite{Mu}, Theorem
2.14), where it is proved for stably free modules.

\begin{theorem}\label{m2}
Let $R$ be a regular affine algebra of dimension $n-1
\geq 2$ over an algebraically closed field $k$. 
Let $A=R[X,f/g]$, where
$g,f$ is a $R[X]$-regular sequence. Suppose 

(1) {\rm char} $k=0$ or {\rm char} $k=p \geq n$.

(2) either $g$ is a monic polynomial or $g(0)\in R^*$. \\
Let $P'$ be a projective $A$-module of rank $n-1$ which is extended
from $R$.  Let $(a,p)\in \Um(A\op P')$ and $P=A\op
P'/(a,p)A$. Then $P \iso P'$.
\end{theorem}

\begin{proof} 
By (\ref{main2}), there exists a projective
$R[X]$-module $Q'$ of rank $n-1$ such that $P=Q'\ot A$. By 
(\ref{lindel1}), $Q'=Q\ot R[X]$ with $Q$ a projective $R$-module
of rank $n-1$. Hence $P=Q\ot_R A$. From (\cite{Mu}, Theorem 2.14), we have
that $K_0(R)\ra K_0(A)$ is injective.  Since $P'$ is extended from $R$ and
$P$ is stably isomorphic to $P'$, hence $Q$ is stably isomorphic to $P'$
as $R$-modules. By (\ref{007}), $Q\iso P'$ as
$R$-modules and hence $P\iso P'$. This proves the result.  $\hfill
\square$ 
\end{proof}


\section{Main Theorem 2}

In this section we prove our second result (\ref{main4}).
Given an affine algebra $A$ over $\BR$
and a subset $I\subset A$, we denote by $Z(I)$, the closed subset of $X
= \Spec A$ defined by $I$ and by $Z_{\tr}(I)$, the set $Z(I)\cap X(\BR)$,
where $X(\BR)$ is the set of 
all real maximal ideals $\mm$ of $A$ (i.e. $A/\mm \iso \BR$).

We begin by stating the following result of Ojanguren and Parimala
(\cite{O-P}, Lemma 2).

\begin{lemma}\label{singular}
Let $A$ be a reduced affine algebra of dimension $n$ over $\BR$ and  
$X=\Spec A$. Let
$[a_1,\ldots,a_d]\in \Um_d(A)$. Suppose $a_1 > 0$ 
on $X(\BR)$. Then, there exists $b_2,\ldots,b_d \in A$
such that $\wt a=a_1+b_2a_2+\ldots+b_da_d > 0$ on $X(\BR)$
and $Z(\wt a)$ is smooth on $X\backslash \Sing X$ of
dimension $\leq n-1$.
\end{lemma}

The following result is analogous to (\cite{O-P},  Proposition 1) and
(\cite{Kesa}, Lemma 3.8).

\begin{lemma}\label{a4}
Let $R$ be a reduced affine algebra of dimension $n\geq 3$ over $\BR$
and let $g,f_1,\ldots,f_r\in R$ with $g$ not belonging to any real
maximal ideal of $R$.
Let $A=R[f_1/g,\ldots,f_r/g]$ and $X=\Spec A$.
Let $P$ be a projective  
$A$-module and let $(a,p)\in \Um(A \op P)$ 
with $a-1\in sgA$ for some $s\in R$.
Then, there exists $h\in 1+gR$ and $\Delta \in
\Aut(A_h \op P_h)$ such that if 
$(a,p)\Delta = (\wt a,\wt p)$, then 

$(1)$ $\wt a > 0$ on $X(\BR)\cap \Spec A_h$,

$(2)$ $Z(\wt a)$ is smooth on $\Spec A_h \backslash \Sing X$ of dimension
$\leq
n-1$ and

$(3)$ $(\wt a,\wt p)=(1,0)$ (mod $ sgA_h$).
\end{lemma}

\begin{proof}
By replacing $g$ by $g^2$, we may assume that $g> 0$ on $X(\BR)$. 
Since $a=1$ mod $sgA$, $(a,sp)\in \Um(A\op P)$. 
Therefore, $a$ has no zero on $Z_{\tr}(O(sp))$. 
Let $r$ be a positive integer such that
$g^ra\in gR$. Let $Y=\Spec R$. Then $g^ra$ has no zero on
$Z_{\tr}(O(sp)) \cap Y(\BR)$.
By {\L}ojasiewicz's inequality ($\cite{BCR}$, 
Proposition 2.6.2), there exists $c\in R$ with $c>0$ on $Y(\BR)$ such
that $1/|a|g^r< c$ on $Z_{\tr}(O(sp)) \cap Y(\BR)$. 
Let $(1+ag^rc)a = a'$. Then $g^ra' >0$ 
on $Z_{\tr}(O(sp)) \cap Y(\BR)$ and hence $a'>0$ on $Z_{\tr}(O(sp))$. Write
$h=1+ag^rc \in 1+gR$. Then $a'=ha$.

Let $W$ be the closed semi-algebraic subset of $X(\BR)$ defined by
$a'\leq 0$. Since
$Z_{\tr}(O(sp)) \cap W =\varnothing$, if $O(p)=(b_1,\ldots,b_d)$ then
$s^2({b_1}^2+\ldots+{b_d}^2) > 0$ on $W$. Hence,
by {\L}ojasiewicz's inequality,
there exists  $c_1\in A$ with $c_1>0$ on $X(\BR)$ such that $|a'|/
g s^2({b_1}^2+\ldots,{b_d}^2) < c_1$. Hence
$a''=a'+c_1g s^2({b_1}^2+\ldots,{b_d}^2)>0$  on $W$ and
hence $a'' >0$ on $X(\BR)$.

We still have $a''=1$ mod $sgA_h$. Since
 $[a'',gs^2{b_1}^2,\ldots,gs^2{b_d}^2] \in \Um_{d+1}(A_h)$,
by (\ref{singular}), there exists 
$h_i\in A_h$ such
 that $\wt a =a'' + \sum _{i=1}^d gs^2{b_i}^2h_i >0$ 
 on $X(\BR)\cap \Spec A_h$ and $Z(\wt a)$ is smooth on 
$\Spec A_h \backslash \Sing X$ of dimension $\leq
n-1$. It is clear from the
proof that there exists $\Delta_1 \in \Aut (A_h\op P_h)$ such that
 $(a,p)\Delta_1 = (\wt a, p)$. Since $\wt a=1$ mod $sgA_h$, there exists
 $\Delta_2 \in E(A_h\op P_h)$ such that $(\wt a,p)\Delta_2 = (\wt
 a,\wt p)$ with $\wt p\in sg P_h$. Take $\Delta=\Delta_1\Delta_2$.
This proves the result. 
 $\hfill \square$
\end{proof}

\begin{lemma}\label{main3}
Let $R$ be an affine algebra of dimension $n\geq 3$ over $\BR$. Let
$g,f_1,\ldots,f_r\in R$ with $g$ a non-zero-divisor and $A=R[f_1/g,\ldots,
f_r/g]$. Assume that
$g$ does not belong to any real maximal ideal of $R$. Let $S=1+gR$ and
$B=A_S$.  Let $P$ be a projective $B$-module of rank $\geq n-1$ which is
extended from $R_S$. Let $(a,p) \in \Um(B\op P)$ with $(a,p)= (1,0)$
mod $Bg$.  Then, there exists  $\wt \Delta \in \Aut(B \op P)$ such
that $(a,p)\wt \Delta = (1,0)$.
\end{lemma}

\begin{proof}
In view of (\ref{bass1}), it is enough to prove the result
when rank of $P$ is $\leq n$. For the sake of simplicity, we assume
that rank of $P=n-1$. The same proof goes through when rank $P=n$.

Without loss of generality, we may assume that $R$ is reduced. 
Let $J_1$ be the ideal of $R_{g}$ defining the
singular locus $\Sing R_{g}$.
Since $R_{g}$ is reduced, $\hh J_1 \geq 1$. Note that $\sqrt
{J_1}=J_1$. Let $J=J_1\cap R$. Then, we may assume that $g$ does not
belong to any minimal primes of $J$ and $\hh J\geq 1$.
Hence $\hh(g,J)R\geq 2$.

Since $(a,p)=(1,0)$ mod $Bg$, by (\ref{a11}), there exists some
$s\in J$ and $\Phi\in \Aut(B\op P)$ such that $P_s$ is free and 
$(a,p)\Phi=(1,0)$ mod $Bsg$. Hence, replacing $(a,p)$ by $(a,p)\Phi$,
we can assume that $(a,p)=(1,0)$ mod $Bsg$. 

There exists some $h\in S$ such that $P$ is a projective $A_h$-module
with $P_s$ free and $(a,p)\in \Um(A_h\op P)$ with $(a,p)=(1,0)$ mod
$sgA_h$. Applying (\ref{a4}),
there exists some $h'\in 1+gR_h$ and $\Delta \in \Aut(A_{hh'} \op
P_{hh'})$ such that $(a,p)\Delta = (a',p')$ with

$(1')~ a'>0$ on $X(\BR)\cap \Spec A_{hh'}$, where $X=\Spec A_h$,

$(2')~ (a',p')=(1,0)$ mod $sgA_{hh'}$ and

$(3')~ Z(a')$ is smooth (since $a'=1$ mod $sgA_{hh'}$ and $s\in
 J_1$) on $\Spec A_{hh'}$ of dimension
$\leq n-1 $.
Note that, since $h^rh' \in 1+Rg$ for some positive integer $r$,
 $A_{hh'}\subset B$. Hence, 
replacing $A_{hh'}$ by $A$ and $(a',p')$ by $(a,p)$,
we assume that the
above properties $(1')-(3')$ holds for $(a,p)$ in the ring $A$, i.e. we
have 

$(1)~ a>0$ on $X(\BR)$, where $X=\Spec A$,

$(2)~ (a,p)=(1,0)$ mod $sgA$ and

$(3)~ Z(a)$ is smooth on $\Spec A$ of dimension
$\leq n-1 $. \\

Since $P_s$ is free of rank $n-1$, there exists a free
submodule $F=A^{n-1}$ of $P$ such that $F_s=P_s$. Replacing $s$ by
a suitable power of $s$, we may assume that $sP\subset F$. Let
$e_1,\ldots,e_{n-1}$ denote the standard basis of $A^{n-1}$.

Since $p\in sgP \subset gF$, $p=b_1e_1+\ldots+b_{n-1}e_{n-1}$ for some
$b_i\in gA$.
Then $[a,b_1,\ldots,b_{n-1}] \in \Um_n(A)$. Let
$T=1+g\BR[g]$ and $C=A_T$. Note that $B=A_S=C\ot C_S$.
Let ``bar'' denotes reduction modulo $Ca$. Since $a-1\in Csg$ and
$s\in J$, $\ol C$
is a smooth affine algebra over $\BR(g)$ of dimension $n-2$. Since
$P_s=F_s$ is free, the inclusion $F\subset
P$ gives rise to equality $\ol F=\ol P$.
In particular, $\ol P$ is free of rank
$n-1 \geq 2$ with a basis $\ol e_1,\ldots,\ol e_{n-1}$ and $\ol p \in
\Um(\ol P)$. 

Assume $n\geq 4$.
We have $[\ol b_1,\ldots,\ol b_{n-1}] \in \Um_{n-1}(\ol
C)$. As in
(\cite{Mu}, Lemma 2.6), by Swan's Bertini theorem (\cite{Swan},
Theorem 1.3), there exists an $\Theta\in E_{n-1}(\ol C)$ such that 
$[\ol b_1,\ldots,\ol b_{n-1}]\Theta=[\ol b_1,\ol b_2,\ol c_3,\ldots,\ol
c_{n-1}]$ with the following properties:

$(1)~ \ol C/\ol J$ is smooth affine $\BR(g)$-algebra
of dimension $2$, where $J$ denotes the ideal of $C$ generated by
$(c_4,\ldots,c_{n-1})$.

$(2)~ \ol C/\ol I$ is smooth affine $\BR(g)$-algebra
of dimension $1$,  where  
$I$ denotes the ideal of $C$ generated by 
$(c_3,\ldots,c_{n-1})$.\\

Every maximal ideal $\mm$
of $\ol C/\ol I$ is the image in $\Spec \ol C/\ol I$
of a prime ideal $\p$ of $C$ of height $n-1$
containing $a$. Since $a$ does not belongs to any real maximal
ideal of $C$, by Serre's result 
\cite{Serre}, the residue field $\BR(\p)=k(\mm)$ of $\mm$  
has cohomological dimension $\leq 1$.
By (\ref{op}), $SK_1(\ol C/\ol I) $ is divisible and the
natural map $K_1Sp(\ol C/\ol I) \ra SK_1(\ol C/\ol I)$ is an isomorphism. 

Let ``tilde'' denotes reduction modulo $\ol I$. Write $D=\ol C$, $\wt
D=D/\ol I$. Then,
there exists $\Theta' \in
\SL_2(\wt D)\cap ESp(\wt D)$ and $t_1,t_2\in D$ such that
$[\wt b_1,\wt b_2] \Theta' = [\wt {t_1}^{n-1},\wt t_2]$.
Since $B=C_S$, $\ol B={\ol C}_S$. We have $\Theta' \in
\SL_2(\ol B/\ol I)\cap ESp(\ol B/\ol I)$ and $t_1,t_2\in \ol B$ such that
$[\wt b_1,\wt b_2] \Theta' = [\wt {t_1}^{n-1},\wt t_2]$.\\

If $n=3$, then $\ol I=0$ and hence $\ol B/\ol I=\ol
B=B/Ba$.
We have $\dim B/Bg=\dim A/Ag=2$ and $B_g$ is an $\BR(g)$-algebra of dimension
two. Thus $\Spec B=\Spec B/Bg \cup \Spec B_g$ with $\dim B/Bg=2=\dim
B_g$. Hence psr$(B)\leq 3$.
Therefore, by (\ref{two}),
 $\Theta'$ has a lift $\Theta_1\in \SL(P\ot B)$.

For $n\geq 4$. Since $\dim \ol B/\ol I \leq 1$ and $\dim \ol B/\ol J
\leq 2$, by (\ref{b1}),
 there exists $\Theta''\in E_{n-1}(B/Ba)$ such
that $[\ol b_1,\ol b_2,\ol c_3,\ldots,\ol c_{n-1}] \Theta'' = [{\ol
t_1}^{n-1},\ol t_2,\ol
c_3,\ldots,\ol c_{n-1}]$. Recall that, there exists $\Theta\in
E_{n-1}(B/Ba)$ such that $[\ol b_1,\ldots,\ol b_{n-1}]\Theta=[\ol
b_1,\ol b_2,\ol c_3,\ldots,\ol c_{n-1}]$.
Since $\ol P$ is free of rank
$n-1\geq 3$,  
$E_{n-1}(\ol A)=E(\ol P)$. By (\ref{one}), 
$\Theta \Theta''\in E_{n-1}(B/Ba)$ can
be lifted to an element $\Theta_1\in SL(P\ot B)$. (In particular, the above
argument shows that every stably free $B/Ba$-module of rank $\geq n-2$ is
cancellative). 

Write $P$ for $P\ot B$.
Thus, in either case ($n\geq 3$), there exists $q \in P$ such that
$$\Theta_1(p)=p_1- aq, ~~{\rm where}~~
p_1={t_1}^{n-1}e_1+t_2e_2+c_3e_3+\ldots+
c_{n-1}e_{n-1}.$$ 

The automorphism $\Theta_1$ of $P$
induces an automorphism $\gL_1=(Id_B,\Theta_1)$ of $B\op
P$. Let $\gL_2$ be the transvection $\Delta_{q}$ of $B\op P$. Then
$(a,p)\gL_1 \gL_2 =(a,p_1)$.

By (\ref{b5}), there exists $\gL_3\in \Aut(B\op P)$ such that
$(a,p_1)\gL_3=(1,0)$. Let $\wt \Delta=\gL_1\gL_2\gL_3$. Then
$\wt \Delta\in
\Aut(B\op P)$ and $(a,p)\wt \Delta=(1,0)$. This proves the result.
$\hfill \square$
\end{proof}

Now, we prove the main theorem of this section.

\begin{theorem}\label{main4}
Let $R$ be an affine algebra of dimension $n\geq 3$
over $\BR$. Let $g,f_1,\ldots,f_r$ be a $R$-regular sequence and
$A=R[f_1/g,\ldots,f_r/g]$. Assume that $g$ does not belong to
any real maximal ideal of $R$. Let $P'$ be a projective $A$-module of rank
$\geq n-1$ which is extended from $R$. Let $(a,p)\in \Um(A\op P')$ 
and $P=A\op P'/(a,p)A$. 
Then, $P$ is extended from $R$.
\end{theorem}

\begin{proof}
By (\ref{a1}), there exists $\Delta\in 
\Aut(A\op P')$ such that $(a,p)\Delta = (1,0)$ mod $Ag$. Let $S=1+Rg$
and $B=A_S$. Applying
(\ref{main3}), there exists $\Delta_1\in \Aut(B\op (P'\ot B))$
such that $(a,p)\Delta\Delta_1 =(1,0)$.
Let $\Psi =
\Delta \Delta_1$. Then, there exists $h\in 1+Rg$ such that $\Psi\in
\Aut(A_h \op {P'}_h)$ and $(a,p)\Psi=(1,0)$.  Rest of the
argument is same as in (\ref{main2}).
$\hfill \square$
\end{proof}

\begin{remark}
The proof of (\ref{main4}) works for any real closed field $k$. For
simplicity, we have taken $k=\BR$.
\end{remark}



\begin{corollary}\label{a3}
Let $R=\BR[X_1,\ldots,X_n]$ and $f,g\in R$ with $g$ not belonging
to any real maximal ideal.  Then, every stably free
$R[f/g]$-modules $P$ of rank $\geq n-1$ is free.
\end{corollary}

\begin{proof}
Write $A=R[f/g]$.
We may assume that $f,g$ have no common factors so that $g,f$ is a regular
sequence in $R$. Since rank $P\geq n-1$, $P\op A^2$ is free. Applying   
(\ref{main4}), we get that $P\op A$ is extended from $R$. By
Quillen-Suslin theorem \cite{Q,Su}, every projective $R$-module is free.
Hence $P\op A$ is free. Again, by (\ref{main4}), $P$ is extended from $R$ and
hence is free.
$\hfill \square$
\end{proof}

The proof of the following result is similar to (\ref{m2}), hence we omit it.

\begin{theorem}
Let $R$ be a regular affine algebra of dimension $n-1 \geq 2$
over $\BR$. Let $A=R[X,f/g]$, where $g,f$ is a $R[X]$-regular sequence.
Suppose  that 

(1) $g$ does not belongs to any real maximal ideal. 

(2) $g$ is a monic polynomial or $g(0)\in R^*$. 
\\
Let $P'$ be a projective $A$-module of rank $n$ which is extended
from $R$.  Let $(a,p)\in \Um(A\op P')$ and $P=A\op
P'/(a,p)A$. Then $P\iso P'$.
 
In particular, every stably free  $A$-module of rank $n$
is free.
\end{theorem}

{\small
{}
}

\end{document}